\documentclass[reqno]{amsart}
\usepackage[pxpazo, numberbasedon=subsection]{zzzams}
\usepackage{bbm}

\DeclareMathOperator{\Bh}{\mathbbm{h}}
\DeclareMathOperator{\cFL}{\mathcal{F}\ell}
\DeclareMathOperator{\Sat}{Sat}
\DeclareMathOperator{\Irr}{Irr}
\DeclareMathOperator{\RTmin}{RTmin}
\DeclareMathOperator{\RT}{RT}
\DeclareMathOperator{\KL}{KL}

\newcommand{\ubar}[1]{\underline{#1}}
\DeclareMathOperator{\cl}{\mathbf{cl}}
\DeclareMathOperator{\AV}{AV}

\DeclareMathOperator{\Rep}{Rep}

\DeclareMathOperator{\trv}{trv}

\newcommand{\laurent}[1]{(\!({\ensuremath{#1}})\!)}
\newcommand{\series}[1]{[\![{\ensuremath{#1}}]\!]}

\title{Affine vertex algebras and an affine analog of Barbasch-Vogan's construction}

\author{Peng Shan$^{1,2}$}
\email{pengshan@tsinghua.edu.cn}

\author{Wenbin Yan$^1$}
\email{wbyan@tsinghua.edu.cn}

\author{Qixian Zhao$^1$}
\email{zhao\_qixian@tsinghua.edu.cn}

\address{\scriptsize{$^1$} Yau Mathematical Sciences Center, Tsinghua University, Beijing, 100084, China}
\address{\scriptsize{$^2$} Department of Mathematical Sciences, Tsinghua University, Beijing, 100084, China}

\begin{document}
\maketitle

This is an expository paper based on the authors' joint work \cite{SYZ:cl} and a few ongoing ones that aim to propose a Langlands style framework for the representation theory of simple affine vertex algebras $L_k(\fg)$ attached to a simple Lie algebra $\fg$ over $\BC$ of simply-laced type and an rational level $k$ greater than the critical level. In this paper, we focus on two conjectures describing the associated variety and simple modules of $L_k(\fg)$ in the case where $k$ is an integer. Our proposal is motivated by classical results on primitive ideals and also by recent new developments in conformal field theory. To convey the ideas, we review the finite dimensional picture and physics background in \S \ref{sec:fin}-\ref{sec:4d} before discussing the results and conjectures in \S \ref{sec:cl}-\ref{sec:conj}.

\subsection*{Acknowledgements}

This work is supported by National Key R\&D Program of China (No. 2025YFA1017400). PS is also supported by NSFC Grant 12225108 and by the New Cornerstone Science Foundation through the Xplorer Prize. The third author would like to thank the \textit{International Congress of Chinese Mathematicians 2025} for the invitation to given a presentation.

\section{Affine vertex algebras}\label{sec:VA}

We let $\fg$ be a simple Lie algebra over $\BC$ with a Cartan subalgebra $\fh$ and a Borel $\fb \supseteq \fh$. We let $\Bh$ (resp. $\check \Bh$) be the Coxeter number (resp. dual Coxeter number) of $\fg$. Write $\rho$ for the half sum of positive roots with respect to the Borel $\fb$. 
The \textit{affine Lie algebra} is a central extension
\begin{equation*}
	0 \aro \BC K \aro \fg_{aff} \aro \fg\laurent{t} \aro 0.
\end{equation*}
Such extensions are determined by a nondegenerate invariant symmetric bilinear form $B(-,-)$ on $\fg$, which we fix to be the unique one such that the induced form on $\fh^*$ satisfies $B(\theta, \theta) = 2$ for the highest root $\theta$. We fix the Iwahori subalgebra to be the preimage of $\fb$ under the map $\fg\series{t} \surj \fg$ sending $t \mapsto 0$. We fix the Cartan subalgebra $\fh_{aff} = \fh \oplus \BC K$ and write its dual as $\fh_{aff}^* = \fh^* \oplus \BC \Lambda_0$ where $\Lambda_0$ is dual to $K$. It contains the weight $\hat \rho := \check \Bh \Lambda_0 + \rho$.

We write $W_{aff}$ for the affine Weyl group, with a set of simple reflections determined by the Iwahori $\fI$. It acts naturally on $\fh_{aff}^*$. The $\circ$-action of $W_{aff}$ on $\fh_{aff}^*$ is given by $w \circ \Lambda = w\cdot(\Lambda+\hat \rho) - \hat \rho$.

\subsection{Associated varieties}

We refer readers to \cite{Kac:VA, Frenkel-Ben-Zvi, Arakawa:W-alg-notes} for the definition of vertex algebras and their properties. Given a vertex algebra $V$ and $a \in V$, we write $Y(a,z) = \sum_{n \in \BZ} a_{(n)} z^{-n-1} \in \End(V) \series{z,z\inv}$ for the corresponding field.

For a number $k \in \BC$, let
\begin{equation*}
	V^k = V^k(\fg) = \cU(\fg_{aff}) \dotimes_{\cU(\fg\series{t}) \oplus \,\BC K} \BC_k,
\end{equation*}
where $\fg\series{t}$ acts trivially on $\BC_k$ and $K$ acts on $\BC_k$ by $k$. This is a highest weight module of $\fg_{aff}$ with highest weight $k \Lambda_0$ and it has the structure of a vertex algebra (see, for example, \cite{Frenkel-Zhu} or \cite[\textsection 2.4]{Frenkel-Ben-Zvi}). 
The unique irreducible quotient $L(k \Lambda_0)$ of $V^k$ as a $\fg_{aff}$-module inherits a structure of a vertex algebra from the one on $V^k$, and the resulting vertex algebra is denoted by
\begin{equation*}
	L_k = L_k(\fg).
\end{equation*}

For a vertex algebra $V$, \textit{Zhu's $C_2$-algebra of $V$} is the ring
\begin{equation*}
	R_V = V/ \{a_{(-2)} b \mid a,b \in V\}
\end{equation*}
with multiplication given by $a \cdot b = a_{(-1)} b$ and a Poisson bracket given by $\{a,b\} = a_{(0)}b$ \cite{Zhu:thesis}. The subspace $\{a_{(-2)} b \mid a,b \in V\}$ is usually denoted by $C_2(V)$, and it coincides with the first piece $F^1 V$ in the filtration on $V$ defined by Li \cite{Li:abVA}. The \textbf{associated variety} of $V$ is the Poisson variety
\begin{equation*}
	X_V := (\Spec R_V )_{red}.
\end{equation*}
We say $V$ is \textit{quasi-lisse} if the Poisson variety $X_V$ has finitely many symplectic leaves \cite{Arakawa-Kawasetsu}. This family of vertex algebras is remarkable, as it has only finitely many simple ordinary modules, and the normalized characters of these modules enjoy modular invariance properties.

In the case of affine vertex algebras, we have $X_{V^k} = \fg^*$ which can (and will) be identified with $\fg$ using the form $B(-,-)$ fixed above. Then $X_{L_k}$ is a closed $G$-stable conic Poisson subvariety inside $X_{V^k} = \fg$, and $L_k$ is quasi-lisse if and only if $X_{L_k} \subseteq \cN$. It is known that $X_{L_k} = \fg$ if $k$ is irrational or if $k \in \BQ$ and $k < -\check \Bh$. When $k = - \check \Bh$, $X_{L_k}$ is equal to the nilpotent cone $\cN$. When $k$ is \textit{admissible} (when $\fg$ is simply-laced, this means $k + \check \Bh$ is a rational number with numerator $\ge \check \Bh$) $X_{L_k}$ is the closure of a single nilpotent orbit. The behavior of $X_{L_k}$ at non-admissible rational levels is interesting but mysterious, and is known only in small families, see \cite{Arakawa-Moreau:Omin, Arakawa-Moreau:sheets, Arakawa-Moreau:irred, AFK, Jiang-Song, ADFLM} for some results in the case where $k$ is integral.

\begin{problem}\label{prob:AV}
	Find a general description of $X_{L_k}$ for any $k \in \BQ_{> - \check \Bh}$.
\end{problem}

In Conjecture \ref{conj:AV} below we propose a uniform description for $X_{L_k}$ for any \textit{integer level} $k > - \check \Bh$ and for all $\fg$ of simply-laced type. For non-integer levels, we refer readers to \cite{SYZ:rat}.

\subsection{Category $\cO$}\label{subsec:O(gaff)}

Consider the category $\cO$ of $\fg_{aff}$ with respect to the Cartan subalgebra $\fh_{aff}$ and the Iwahori subalgebra $\fI \subset \fg_{aff}$, see \cite[\S 9.1]{Kac:inf-dim-alg}. For an integral weight $\Lambda \in \fh_{aff}^*$, write $L(\Lambda)$ for the irreducible object in $\cO(\fg_{aff})$ with highest weight $\Lambda$, write $\cO_\Lambda(\fg_{aff})$ for the block whose simple objects are of the form $L(y \circ \Lambda)$ for $y \in W_{aff}$, and write $K_0 \cO_\Lambda(\fg_{aff})$ for its Grothendieck group. If $\Lambda+ \hat \rho$ is dominant but possibly singular, that is $\langle \check \alpha_i, \Lambda + \hat \rho \rangle \in \BZ_{\ge 0}$ for any affine simple coroot $\check \alpha_i$, then we have a bijection
\begin{align*}
	\Irr \cO_\Lambda(\fg_{aff}) \bij& \{y \in W_{aff} \mid y \text{ is the longest element in } y W_\Lambda \},\\
	L(y \circ \Lambda) \mathop{\mapsfrom}& y
\end{align*}
where $W_\Lambda$ is the stabilizer of $\Lambda$ in $W_{aff}$ under the $\circ$-action, see, for example, \cite{Deodhar-Gabber-Kac,KT:KL2,KT:KL3}.

It is known that the category of $V^k$-modules is a full subcategory of smooth $\fg_{aff}$-modules \cite[Lemma 4.3.2, Theorem 5.1.6]{Frenkel-Ben-Zvi}, and the essential image consists of those modules on which the central element $K \in \fg_{aff}$ acts by $k$. Here, a $\fg_{aff}$-module $M$ is \textit{smooth} if for any $v \in M$ and any $a \in \fg$, there exists an $N > 0$ so that $a t^n \cdot v = 0$ for all $n > N$. In particular, any object in the category $\cO$ is smooth, and we have $\cO(V^k) \subseteq \cO(\fg_{aff})$. Via the quotient map $V^k \surj L_k$, the category of $L_k$-modules is naturally a full subcategory of $V^k$-modules. We thus obtain the category $\cO(L_k)$ as the intersection of $\cO(V^k)$ with the category of $L_k$-modules. We similarly have the subcategories
\begin{equation*}
	\cO_\Lambda(L_k) \subseteq \cO_\Lambda(V^k) \cong \cO_\Lambda(\fg_{aff})
\end{equation*}
for any element $\Lambda \in k \Lambda_0 + \fh^* \subset \fh_{aff}^*$.

\begin{problem}\label{prob:O}
	Find a general description of $\Irr \cO(L_k)$ as a subset of $W_{aff}$.
\end{problem}

Existing approaches to this problem often go through careful calculations with primitive vectors in $V^k$ or the Zhu's algebra $A(L_k)$ \cite{Zhu:modular}, which are difficult to generalize. In Conjecture \ref{conj:O(Lk)} below we propose a description of simple modules in the cases where $L_k$ is conjecturally quasi-lisse. It is stated in terms of something (namely Kazhdan-Lusztig left cells) which, we think, are more conceptual.

\section{The classical story}\label{sec:fin}

In this section we explain the classical prototype of our constructions. The main takeaway is the diagram (\ref{diag:fin-duality}). Let us write $W$ for the (finite) Weyl group. The dot-action $W \acts_\circ \fh^*$ is given by $w \circ \lambda = w (\lambda+\rho) - \rho$. Via the Harish-Chandra isomorphism, infinitesimal characters $\chi_\lambda: \cZ(\cU(\fg)) \to \BC$ are parameterized by $(W,\circ)$-orbits $W \circ \lambda$ in $\fh^*$. We follow the convention that $\chi_0$ is the infinitesimal character of the trivial representation.

\subsection{Kazhdan-Lusztig cells}\label{subsec:KL-cells}


In this subsection we let $(W,S)$ denote a Coxeter group. Later $W$ will become the finite Weyl group (again denoted by $W$) or the affine Weyl group $W_{aff}$. Let $\cH$ denote the Hecke algebra of $W$ over $\BZ[q^{\pm1/2}]$, a free $\BZ[q^{\pm1/2}]$-module with basis $T_w$, $w \in W$ satisfying
\begin{equation*}
	T_x T_y = T_{xy} \text{ if } \ell(x) + \ell(y) = \ell(xy),
\end{equation*}
\begin{equation*}
	(T_s + 1)(T_s - q) = 0 \text{ for any simple reflection } s.
\end{equation*}
There is an involution $h \mapsto \bar h$ of $\cH$ defined by
\begin{equation*}
	\bar q^{1/2} = q^{-1/2}, \quad \bar T_w = T_{w\inv}\inv.
\end{equation*}
The Kazhdan-Lusztig basis elements $C_w$ of $\cH$ are given by
\begin{equation*}
	C_w = q^{-\ell(w)/2} \sum_{y \le w} (-q)^{\ell(w) - \ell(y)} \bar P_{y,w} T_y,
\end{equation*}
for some uniquely determined $P_{y,w} \in \BZ[q]$ with degree $\le \frac12 (\ell(w) - \ell(y) - 1)$ satisfying $\bar C_w = C_w$. See \cite[Theorem 1.1]{KL:Hecke}. 

Whenever $s \in W$ is a simple reflection, $w \in W$, and $sw > w$, the product $C_s C_w$ is a $\BZ$-linear combination of $C_y$'s for some $y \le sw$. If $C_y$ appears, we say $y \le_L w$. We use $\le_L$ to generate a preorder on $W$. If $y \le_L w$ and $y \ge_L w$, we write $y \sim_L w$. The equivalence classes of $\sim_L$ are called \textit{left cells}
\begin{equation*}
	\bc^L(w) = \{ y \in W \mid y \sim_L w\}.
\end{equation*}
The same definition can be made using $C_w C_s$ ($ws > w$) instead, the resulting relations are denoted $\le_R$, $\sim_R$, and the equivalence classes are called \textit{right cells}, denoted by $\bc^R(w)$. The preorder generated by using both $\le_L$ and $\le_R$ is denoted $\le_{LR}$, and we have the corresponding equivalence relation $\sim_{LR}$ and the equivalence classes $\ubar \bc(w)$, the \textit{two-sided cells}. The relation $\le_R$ induces a filtration on $\cH$ by right ideals
\begin{equation*}
	\cH_{\le_R w} = \bigoplus_{y \le_R w} \BZ[q^{\pm1/2}] \cdot C_y \subset \cH.
\end{equation*}
The subquotients are denoted by $\cH_{\bc^R(w)}$ and are called the \textit{right cell modules}. After specializing at $q=1$, they become modules over $\cH|_{q=1} = \BZ[W]$, i.e. representations of $W$.

There is a ring involution on $\cH$ given by
\begin{equation*}
	{}^*(-) : \cH \to \cH,\quad
	{}^*q^{1/2} = - q^{1/2},\quad
	{}^*T_w = (-q)^{\ell(w)} T_{w\inv}\inv.
\end{equation*}
Let $\cH^* = \Hom_{\BZ[q^{\pm1/2}]}(\cH,\BZ[q^{\pm1/2}])$. Let $D_w \in \cH^\vee$ be defined by
\begin{equation*}
	\langle D_w, C_y \rangle = \delta_{w,y\inv},
\end{equation*}
where $\langle -,- \rangle$ is the natural pairing between $\cH^*$ and $\cH$.
There is a left $\cH$-module structure on $\cH^*$ given by
\begin{equation}\label{eqn:H-acts-on-H*}
	\langle h \cdot f, h' \rangle = \langle f, h' \cdot {}^*h \rangle \text{ for any } f \in \cH^\vee, h,h' \in \cH.\footnotemark
\end{equation}
\footnotetext{One may identify $\cH^*$ with $\hat \cH := \prod_{y \in W} \BZ[q^{\pm1/2}] \cdot T_y$, and $\langle h_1,h_2 \rangle$ with $\tau( h_1 h_2)$, where $\tau: \hat \cH \to \BZ[q^{\pm1/2}]$ sends $T_1 \mapsto 1$ and $T_w \mapsto 0$ for $w \neq 1$. Then $D_w$ becomes the element defined in (say) \cite[1.3]{Lusztig:aff-cells-1}, and the left module structure $\cH \acts \cH^*$ becomes the left multiplication action on $\hat \cH$ composed with the involution on $\cH$ given by $T_w \mapsto (-q)^{\ell(w)} T_{w\inv}\inv$.}%
For any simple reflection $s$, $C_s \cdot D_y$ is a (finite) $\BZ[q^{\pm1/2}]$-linear combination of $D_z$'s satisfying $z \ge_L y$, see \cite[(5.1.16)]{Lusztig:char}. As a result, the subspace
\begin{equation*}
	\cH^\vee := \bigoplus_{y \in W} \BZ[q^{\pm 1/2}] \cdot D_y \subseteq \cH^*
\end{equation*}
is an $\cH$-submodule in $\cH^*$, and we have a filtration on $\cH^\vee$ by left $\cH$-submodules given by
\begin{equation*}
	\cH_{\ge_L w}^\vee = \bigoplus_{y \ge_L w} \BZ[q^{\pm1/2}] \cdot D_y,
\end{equation*}
called the left cell filtration. The subquotients are denoted by $\cH_{\bc^L(w)}^\vee$, called the \textit{dual left cell modules}. We also write
\begin{equation*}
	\cH_{\not\le_L w}^\vee 
	:= 
	\bigoplus_{\substack{%
			y \in W\\
			y \not\le_L w}}
	\BZ[q^{\pm1/2}] \cdot D_y.
\end{equation*}
This is a sum of cell submodules $\cH_{\ge_L \bc^L}^\vee$ over those left cells $\bc^L$ that are either $>_L w$ or not comparable to $w$ under $\le_L$. Then the natural inclusion $\cH_{\ge_L w}^\vee \inj \cH^\vee$ induces an inclusion of $\cH$-modules
\begin{equation}\label{eqn:cell-in-quotient}
	\cH_{\bc^L(w)}^\vee \injects \cH^\vee/ \cH_{\not\le_L w}^\vee;
\end{equation}
the right hand side has a basis given by images of $D_y$, $y \le_L w$.

Finally, the pairing $\langle -,- \rangle$ between $\cH^\vee$ and $\cH$ induces a perfect pairing 
\begin{equation}\label{eqn:pairing}
	\cH_{\bc^L(w)}^\vee \times \cH_{\bc^R(w\inv)} \aro \BZ[q^{\pm1/2}],
\end{equation}
see \cite[Proposition 1.2.2]{Bonnafe}.

\subsection{Springer correspondence}

We write $\cN$ for the nilpotent cone of $\fg$ and $\ubar \cN$ for the set of adjoint orbits. Let $e \in \BO \in\ubar \cN$. The \textit{Springer fiber} of $e$ is 
\begin{equation*}
	\cB_e = \big\{ g B \in G/B \mid \Ad(g)\inv e \in \fb \big\},
\end{equation*}
a closed subvariety in the flag variety $G/B$. Springer constructed in \cite{Springer:W} a $W$-action\footnote{Here we use the convention of Lusztig \cite{Lusztig:sp-1}, which differs from Springer's original definition by tensoring with the sign representation. In particular, the zero orbit (resp. principal orbit) corresponds to the sign representation (resp. trivial representation).} on the cohomologies $H^*(\cB_e)$ that commutes with the natural $A_G(e)$-action, where $A_G(e)$ is the component group of the centralizer $Z_G(e)$ of $e$ in $G$. In particular, for any irreducible representation $\tau$ of $A_G(e)$, or equivalently a $G$-equivariant irreducible local system $\tau$ on the orbit $\BO$, the $\tau$-isotypic part of the top cohomology $H^{top}(\cB_e)^\tau$ is a representation of $W$. It is either irreducible or zero, and they exhaust all irreducible representations of $W$. This leads to an injection, the \textit{Springer correspondence}
\begin{equation*}
	\Irr \Rep(W) \injects \big\{ (\BO, \tau) \mid \BO \in \ubar \cN, \tau \in \Irr A(\BO) \big\}
\end{equation*}
where $A(\BO)$ is $A_G(e)$ for some $e \in \BO$ and $\Irr$ denotes the set of irreducible representations. For $(\BO,\tau)$ in the image, we write $E_{(\BO,\tau)}$ for the corresponding $W$ representation.

In \cite{Lusztig:sp-1} a class of irreducible representations of $W$ is defined, later called \textit{special representations} of $W$. We refer to \textit{loc. cit.} for its definition. Under the Springer correspondence, special representations always correspond to pairs $(\BO, \trv)$ for some orbit $\BO$ and the trivial representation of $A(\BO)$. An orbit $\BO$ is said to be \textit{special} if it arises in this way. We write
\begin{equation*}
	\ubar \cN_{sp} = \big\{ \text{special orbits in } \ubar \cN \big\}.
\end{equation*}


\subsection{Primitive ideals and left cells}

A \textit{primitive ideal} $J$ in $\cU(\fg)$ is the annihilator of a simple $\fg$-module $L$. We say $J = \Ann L$ has infinitesimal character $\chi_\lambda$ if $L$ does. The \textit{associated variety} of $J$, denoted by $\AV J$, is the associated variety of the $\fg$-module $\cU(\fg)/J$, the latter can be defined as the support of the $\Sym(\fg) = \BC[\fg^*]$-module $\gr (\cU(\fg)/J)$, where the associated graded is taken with respect to a good filtration on $\cU(\fg)/J$, see \cite[Appendix D]{HTT}. After identifying $\fg^*$ with $\fg$, $\AV(J)$ is a closed subvariety in $\fg$. 

\begin{theorem}~
	\begin{enumerate}[label=(\roman*)]
		\item (Duflo \cite{Duflo:prim}) Every primitive ideal is the annihilator of a simple highest weight module $L(\mu)$.
		
		\item (Conze–Dixmier \cite{Conze-Dixmier}) There is a unique maximal primitive ideal $J_{\lambda,max}$ for each infinitesimal character $\chi_\lambda$, and $J_{\lambda,max} = \Ann L(\lambda)$ if $\lambda$ is dominant.
		
		\item (Borho-Brylinski \cite{Borho-Brylsinki:diff-op-1}, Joseph \cite{Joseph:prim}, Vogan \cite{Vogan:AV}) The associated variety $\AV J$ of a primitive ideal is the closure of a single nilpotent orbit.
	\end{enumerate}
\end{theorem}

For $\lambda + \rho \in \fh^*$ dominant integral, let $\cO(\fg)_\lambda$ denote the category $\cO$ for $\fg$ with infinitesimal character $\chi_\lambda$, see \cite{Humphreys:Cat-O}. It contains Verma modules $M(w \circ \lambda)$ for $w \in W$ and their irreducible quotients $L(w \circ \lambda)$. Let $W_\lambda$ be the stabilizer of $\lambda$ under the $(W,\circ)$-action. As $w$ ranges through longest coset representatives of $W/W_\lambda$, the Vermas and the simples form two bases of the Grothendieck group $K_0 \cO(\fg)_\lambda$. It is elementary to show that $w$ is longest in $w W_\lambda$ if and only if $w \le_L w_\lambda$ where $w_\lambda \in W_\lambda$ is the longest element. We thus obtain an isomorphism of abelian groups
\begin{equation}\label{eqn:O-Hecke}
	K_0 \cO(\fg)_\lambda \cong \big( \cH^\vee / \cH_{\not\le_L w_\lambda}^\vee \big)|_{q=1},\quad
	L(y \circ \lambda) \mathop{\mapsfrom} D_y|_{q=1}.
\end{equation}
In fact this isomorphism can be made into $W$-linear, where the $W$-action on the left side is constructed either by using twisting functors \cite{Andersen-Stroppel:twisting} or by using Beilinson-Bernstein's localization theorem and the left convolution action of $D^b(B \backslash G/B)$ on itself \cite[Chapter 12-13]{HTT}.

Write $\cO(\cU(\fg)/J)$ for the full subcategory of $\cO(\fg)_\lambda$ consisting of modules annihilated by a primitive ideal $J$ of infinitesimal character $\lambda$. Its Grothendieck group is a $\cH|_{q=1}$-submodule of $K_0 \cO(\fg)_\lambda$. The following fact is the outcome of developments of Springer theory and the classification of primitive ideals and their associated varieties, the full history of which we do not attempt to cover. Relevant references include \cite{Joseph:W-mod-prim,Joseph:prim,Vogan:ordering,Lusztig:char,Barbasch-Vogan:unipotent,KL:Hecke,Beilinson-Bernstein:localization,Brylinski-Kashiwara}.

\begin{theorem}\label{thm:O(J)}
	Suppose $\lambda+ \rho$ is dominant integral. Assume $x,y \le_L w_\lambda$, i.e. $x$ and $y$ are longest in their $W_\lambda$-cosets.
	\begin{enumerate}[label=(\roman*)]
		\item $L(x \circ \lambda)$ and $L(y \circ \lambda)$ have the same annihilator if and only if $x \sim_L y$. The isomorphism (\ref{eqn:O-Hecke}) restricts to
		\begin{equation*}
			K_0 \cO(\cU(\fg)/ \Ann L(x \circ \lambda)) \cong \big( \cH_{\ge_L x}^\vee / \cH_{\ge_L x, \not\le_L w_\lambda}^\vee \big)|_{q=1}
		\end{equation*}
		(the right side is spanned by $D_y|_{q=1}$ with $x \le_L y \le_L w_\lambda$). In particular,
		\begin{equation*}
			K_0 \cO(\cU(\fg)/J_{\lambda,max}) \cong \cH_{\bc^L(w_\lambda)}^\vee|_{q=1}.
		\end{equation*}
		
		\item $\Ann L(x \circ \lambda)$ and $\Ann L(y \circ \lambda)$ have the same associated variety if and only if $x \sim_{LR} y$. For $\lambda+\rho$ regular, the assignment sending $x$ to the unique open nilpotent orbit inside $\AV \Ann L(x \circ \lambda)$ defines an order-reversing bijection
		\begin{equation*}
			\{\text{two-sided cells in } W\} \bijects \ubar \cN_{sp}.
		\end{equation*}
		
		\item For any $w \in W$, there is a unique special representation appearing in the dual left cell module $\cH_{\bc^L(w)}^\vee|_{q=1}$ with multiplicity one. It is given by $E_\BO = E_{(\BO,\trv)}$ where $\BO$ is the image of $\ubar \bc(w)$ under the bijection in (ii).
	\end{enumerate}
\end{theorem}

\subsection{Barbasch-Vogan-Lusztig-Spaltenstein duality}\label{subsec:BV-dual}

We now come to the key classical construction, a duality map between nilpotent orbits due to Barbasch-Vogan \cite{Barbasch-Vogan:unipotent} and earlier using a different method by Lusztig-Spaltenstein \cite{Lusztig:sp-1}.

Let $\check \fg$ be the Langlands dual Lie algebra containing a Cartan $\check \fh$ which is identified with $\fh^*$. We write $\check \cN$ for the nilpotent cone of $\check \fg$ and write $\ubar{\check \cN}$ for the set of orbits in it. 
The duality map
\begin{equation*}
	\bd: \ubar{\check \cN} \aro \ubar{\cN}
\end{equation*}
is defined as follows. For any orbit $\check \BO \in \ubar{\check \cN}$, by the theorem of Jacobson-Morozov there exists an $\fsl_2$-triple $\{h,e,f\}$ with $e \in \check \BO$ and $h \in \check \fh = \fh^*$ dominant. Then the $(W,\circ)$-orbit of $\lambda := \frac h2 - \rho \in \fh^*$ determines an infinitesimal character $\chi_\lambda$ and hence a maximal primitive ideal $J_{\lambda,max}$ in $\cU(\fg)$. Let $\BO_\lambda$ be the unique open orbit in $\AV J_{\lambda,max}$. Then we set
\begin{equation*}
	\bd \check \BO := \BO_\lambda.
\end{equation*}
It is an order-reversing surjection onto the set of special orbits $\ubar \cN_{sp}$. Each fiber of $\bd$ contains a unique maximal orbit, and it is special. As a result, $\bd$ restricts to a bijection $\ubar{\check \cN}_{sp} \bij \ubar \cN_{sp}$.

Putting the map $\bd$ and the bijection in Theorem \ref{thm:O(J)}(ii) together, we obtain a commutative triangle
\begin{equation}\label{diag:fin-tri}
	\begin{tikzcd}[column sep=small]
		& \{\text{two-sided cells in }W\} \ar[dr, leftrightarrow, "{\ref{thm:O(J)}(ii)}"] \\
		\ubar{\check \cN}_{sp} \ar[rr, leftrightarrow, "\bd"] \ar[ur, "(\ast)"]
		&& \ubar \cN_{sp}
	\end{tikzcd}
\end{equation}
If we start with an \textit{even} orbit $\check \BO \in \ubar{\check \cN}_{sp}$, that is if the weight $\lambda = \frac h2 - \rho$ constructed from $\check \BO$ is integral, then by tracing through the definitions one quickly sees that the top-left arrow $(\ast)$ sends 
\begin{equation}\label{eqn:even-to-cell}
	(\ast): \check \BO \mapsto \ubar \bc(w_\lambda) \qquad (\text{for } \check \BO \text{ even}).
\end{equation}
Moreover, the Springer representation $E_{\check\BO}$ is the unique special representation appearing in the right cell module $\cH_{\bc^R(w_\lambda)}|_{q=1}$, and it appears with multiplicity one. In particular, we have a unique surjection of $W$ representations 
\begin{equation}\label{eqn:cR-to-Spr}
	\cH_{\bc^R(w_\lambda)}|_{q=1} \surj E_{\check \BO} = H^{top}(\check \cB_e)^{A_{\check G}(e)}	
\end{equation}
for any $e \in \check \BO$.

\subsection{The classical picture}\label{subsec:classical-dual}

Putting everything together, we obtain the following picture for an even orbit $\check \BO \in \ubar{\check \cN}_{sp}$:

\begin{equation}\label{diag:fin-duality}
	\begin{tikzcd}
		\cH_{\bc^R(w_\lambda)}|_{q=1} \ar[d, two heads, "(\ref{eqn:cR-to-Spr})"'] \ar[rr, leftrightarrow, "\text{perfect pairing}", "(\ref{eqn:pairing})"']
		&& \cH_{\bc^L(w_\lambda)}^\vee|_{q=1} \ar[d, equal, "\ref{thm:O(J)}(i)"]
		\\
		H^{top}(\underbracket{\check \cB_e})^{A_{\check G}(e)} \ar[d, "\text{moment map}"']
		& \ubar \bc(w_\lambda) \ar[dl, leftrightarrow, "\ref{eqn:even-to-cell}", "(\ast)"'] \ar[dr, leftrightarrow, "\ref{thm:O(J)}(ii)"'] \ar[ul, dotted] \ar[ur, dotted]
		& K_0 \cO(\underbracket{\cU(\fg)/J_{\lambda,max}}) \ar[d, "\AV"]
		\\		
		\check \BO \ar[rr, leftrightarrow, "\bd"] 
		&& \BO
	\end{tikzcd}
\end{equation}
Here the arrow ``moment map'' is simply a sloppy way of indicating that the Springer fiber $\check \cB_e$ lies over the orbit $\check \BO$. More precisely, the Springer fiber $\check \cB_e$ can be defined as the fiber of $e \in \check \cN$ of the moment map $\mu:T^* (\check G/ \check B) \surj \check \cN$, and so $\mu(\check \cB_e) = \{e\}$, which lies in the orbit $\check \BO$. The dotted arrows going out of $\ubar \bc(w_\lambda)$ simply mean that the two modules $\cH_{\bc^R(w_\lambda)}|_{q=1}$ and $\cH_{\bc^L(w_\lambda)}^\vee|_{q=1}$ are defined using the right/left cell of $w_\lambda$, both are contained in the two-sided cell $\ubar \bc(w_\lambda)$.

Our main conjectures proposes an affine analog of this diagram, see (\ref{diag:aff-duality}).

\section{Motivation from 4D mirror symmetry}\label{sec:4d}

We discuss briefly some recent developments coming from mathematical physics that also motivated our work. 

Although vertex operator algebras (VOAs) are known to provide a rigorous mathematical framework for two dimensional chiral conformal field theories, a recent breakthrough discovery by Beem-Lemos-Liendo-Peelaers-Rastelli-van Rees \cite{BLLPRvR} showed that they also appear as almost complete invariants of four dimensional $\cN=2$ super conformal field theories (SCFT). More precisely, they construct a map $\BV$ sending a 4D $\cN =2$ SCFT $\cT$ to a VOA $\BV(\cT)$. Conjecturally this map is finite-to-one, and the Higgs branch of $\cT$ is expected to coincide with the associated variety $X_{\BV(\cT)}$ of $\BV(\cT)$ \cite{Beem-Rasterlli:VOA-Higgs}. Moreover, the VOAs arising in this context must be quasi-lisse. Since the simple affine vertex algebras $L_k(\fg)$ sometimes appear as $\BV(\cT)$ or as extensions of $\BV(\cT)$'s, it is essential from this perspective to determine which $L_k(\fg)$ is quasi-lisse.

Mirror symmetries also come into this picture. For 3D supersymmetric gauge theories, mirror symmetry involves two hyperk\"ahler manifolds $X$ and $Y$. By the work of Braden-Licata-Proudfoot-Webster \cite{BLPW:Gale, BLPW:symp-dual-I, BLPW:symp-dual-II}, the relation between $X$ and $Y$ can be partly understood as a duality between geometric information extracted from $Y$ and representation theoretic data from the category $\cO$ of a quantization $\cA_X$ of $X$ (and its resolution). In a similar flavor, Shan-Xie-Yan \cite{Shan-Xie-Yan:phys} proposed a mirror symmetry for 4D $\cN=2$ SCFTs, where the role of $\cA_X$ is played by the VOA $\BV(\cT)$ and the role of $Y$ is played by the Coulomb branch of $\cT$. The latter is the moduli space $\cM_\gamma$ of certain Higgs bundles \cite{BBAMY:wild-Higgs}, which, when considering topological properties, can be replaced by a homogeneous affine Springer fiber associated to the Langlands dual algebra of $\hat \fg$. Many aspects of this proposal were confirmed in the subsequent work \cite{Shan-Xie-Yan} of Shan-Xie-Yan in the case of boundary admissible levels. However, the general case of homogeneous levels, as well as extensions of the proposal to non-homogeneous levels, remain largely unclear.

The objects mentioned here form the following diagram
\begin{equation}\label{diag:4d-mirror}
	\begin{tikzcd}
		&& \text{4d $\cN=2$ SCFT } \cT \ar[dl, "\text{Coulomb branch}"'] \ar[ddr, "\text{Higgs branch (conj.)}", near start, distance=2.5cm, out=5, in=40] \ar[dr, "\BV"]
		\\
		\cFL_\gamma \ar[r, hook] 
		& \cM_\gamma \ar[rr, leftrightarrow, "\text{``dual''}"]  
		&& \text{VOA } \BV(\cT) \ar[d, "\AV"']
		\\
		&&& X_{\BV(\cT)}
	\end{tikzcd}
\end{equation}
and it looks reminiscent to parts of the diagram (\ref{diag:fin-duality}). This motivated us to construct the rest of (\ref{diag:fin-duality}) in this new setting for any integer level above critical.

\section{Cyclotomic levels}\label{sec:cl}

We now turn to the affine analog of (\ref{diag:fin-duality}), and to our proposal for Problems \ref{prob:AV} and \ref{prob:O}. We retain the notations defined in \S\ref{sec:VA}. 

\subsection{\texorpdfstring{The map $\cl_n$ and the orbits $\BO(m)$}{The map cln and the orbits O(m)}}

The following map is the key construction in our method.

\begin{definition}[{\cite[Definition 2.1.2]{SYZ:cl}}]
	The \textbf{(nilpotent) cyclotomic level map}
	\begin{equation*}
		\cl_n: \cN \aro \BZ_{\ge 1}
	\end{equation*}
	is defined as follows. For $e \in \cN$, let $\fl \subset \fg$ be a Bala-Carter Levi of $e$ (a minimal Levi subalgebra containing $e$), and let $\{h,e,f\}$ be an $\fsl_2$-triple in $\fl$. Then the action $\ad h \acts \fl$ has eigenvalues in $2\BZ$. Let $2a$ be the largest eigenvalue of $\ad h \acts \fl$. We then set
	\begin{equation*}
		\cl_n(e) := a+1.
	\end{equation*}	
	This map has been considered before in various specific settings, for example in \cite[Proposition 9.2]{KL:Fl} and in \cite[\S7.3]{stable-grading} for distinguished nilpotent elements. Alternatively, we set
	\begin{equation*}
		\cl_n(e) := \min\{m \in \BZ_{\ge 1} \mid (\ad e)^{2m} = 0\}.
	\end{equation*}
	The equivalence of the two definitions is proven in \cite[Lemma 2.1.4]{SYZ:cl}, using a case-by-case analysis of the largest $\ad h$-eigenvalue on $\fl$ and $\fg$. The map $\cl_n$ clearly descends to a map on orbits
	\begin{equation*}
		\cl_n: \ubar \cN \aro \BZ_{\ge 1}.
	\end{equation*}
	It is order preserving (as can be easily seen from its second definition) and the maximal value is the Coxeter number $\Bh$, achieved by the principal orbit.
\end{definition}

The following remarkable property of $\cl_n$ is taken from \cite[Theorem 2.1.6]{SYZ:cl}. In view of the second definition of $\cl_n$, the theorem follows from \cite{VIGRE}. The proof in \textit{op. cit.}, however, relies again on a case-by-case examination.

\begin{theorem}\label{thm:Om}
	For any positive integer $m$, there exists a unique orbit $\BO(m) \in \ubar \cN$ so that
	\begin{equation*}
		\cl_n\inv([1,m]) = \overline{\BO(m)}.
	\end{equation*}
\end{theorem}

Recall that each fiber of the duality map $\bd$ contains a unique maximal orbit, and it is special. We exhibit a similar behavior here from the uniqueness of $\BO(m)$. The explicit values of $\cl_n$ on each orbit and the description of $\BO(m)$ for each $m$ are contained in \cite[\S2.1]{SYZ:cl}. In classical types, $\cl_n(\BO)$ equals either $p_1$ or $p_1-1$, where $p_1$ is the largest part of the partition $p$ parameterizing $\BO$, and $\BO(m)$ corresponds roughly to the partition that is most rectangular with width $m$.

\subsection{Relation to two-sided cells}\label{subsec:cl-cells}

We next move on to the bijection $(\ast): \ubar{\check \cN}_{sp} \bij \{\text{two-sided cells in }W\}$ mentioned in (\ref{diag:fin-tri}). The notion of Kazhdan-Lusztig cell reviewed in \S\ref{subsec:KL-cells} make sense for $W_{aff}$, and the affine analog of the map $(\ast)$ was constructed by Lusztig \cite[Theorem 4.8]{Lusztig:aff-cells-4}
\begin{equation}\label{eqn:Lusztig-bij}
	\ubar{\check \cN} \bijects \{\text{two-sided cells in } W_{aff}\}
\end{equation}
which subsumes the bijection $(\ast)$: the subset $\ubar{\check \cN}_{sp}$ of special orbits correspond to the set of cells in $W_{aff}$ that meet the finite Weyl group $W$, and taking intersection with $W$ recovers the bijection $(\ast)$. However, the realization (\ref{eqn:even-to-cell}) of the map $(\ast)$ on even orbits does not seem to have a known affine analog in general. Our next theorem can be viewed as an affine analog to (\ref{eqn:even-to-cell}) on orbits of the form $\check \BO(m)$.

For each $m \in \BZ_{\ge 1}$, consider the weight $m \Lambda_0 + \rho \in \fh_{aff}^*$. By an affine rho-shift, $(m \Lambda_0 + \rho) - \hat \rho = k \Lambda_0$ is the highest weight of the simple affine vertex algebra $L_k(\fg) = L(k \Lambda_0)$, viewed as a $\fg_{aff}$-module. Note however that $m \Lambda_0 + \rho$ is in general not dominant. Let
\begin{itemize}
	\item $\xi_m \in \fh_{aff}^*$ be the unique dominant $W_{aff}$-translate of $m \Lambda_0 + \rho$ under the usual (non-$\circ$) action (so that the dominant translate of $k \Lambda_0$ under the $\circ$-action is $\xi_m - \hat \rho$), let 
	
	\item $W_m \subset W_{aff}$ be the stabilizer of $\xi_m$, let
	\item $w_m \in W_m$ be the longest element, and let
	
	\item $\ubar \bc(w_m) \subset W_{aff}$ be the two-sided cell containing $W_{aff}$.
\end{itemize}
The following theorem is proven in \cite[Theorem 4.2.1]{SYZ:cl}.

\begin{theorem}\label{thm:Om-wm}
	Let $m \in \BZ_{\ge 1}$, and let $\fg$ be of simply-laced type. Under Lusztig's bijection (\ref{eqn:Lusztig-bij}), the orbit $\check \BO(m)$ corresponds to the cell $\ubar \bc(w_m)$.
\end{theorem}

The assignment $\check \BO(m) \xmapsto{\cl_n} m \mapsto \xi_m \mapsto \ubar \bc(w_m)$ is very similar to the restriction of $(\ast)$ to even orbits, which sends $\check \BO \mapsto \lambda =\frac h2 - \rho \mapsto \ubar \bc(w_\lambda)$. According to \cite[Proposition 4.8]{Yun:Epi} (see also \cite[Lemma 4.2.3]{SYZ:cl}), the theorem is equivalent to saying that under the Springer correspondence, the pair $(\check \BO(m),\trv)$ is mapped to Lusztig-Spaltenstein's $j$-induction $j_{W_m}^W sgn$ of the sign representation of $W_m$. This is verified on a case-by-case basis. Along the way, we computed explicitly the dominant weights $\xi_m$ for all $1 \le m \le \Bh$, which seem to be of independent interest.

\subsection{Relation to Springer fibers}

Finally, we discuss a relation between the map $\cl_n$ and affine Springer fibers, as well as the reason for the term ``cyclotomic'' in the terminology. Since Springer fibers appear in both the classical picture (\ref{diag:fin-duality}) and in 4D mirror symmetry (\ref{diag:4d-mirror}), the map $\cl_n$ should have some natural compatibility with them. 

Let $\cFL = G\laurent{t}/I$ be the affine flag variety, where $I$ is the Iwahori subgroup. For an element $\gamma \in \fg\laurent{t}$, the \textit{affine Springer fiber} attached to $\gamma$ is 
\begin{equation*}
	\cFL_\gamma = \{gI \in \cFL \mid \Ad(g)\inv \gamma \in \fI^\circ \} \subseteq \cFL
\end{equation*}
where $\fI^\circ$ is the radical of the Iwahori subalgebra $\fI$. We assume $\gamma$ is \textit{regular semisimple} (rs for short), that is, $\gamma$ is regular semisimple when viewed as an element of $\fg(\bar F)$, where $\bar F = \bigcup_{i \ge 0} \BC \laurent{t^{1/i}}$ is the algebraic closure of $F = \BC\laurent{t}$. We also assume $\gamma$ is \textit{topologically nilpotent} (tn for short), that is, $\gamma$ can be conjugated under $G\laurent{t}$ to an element in $\fI^\circ$. Under these assumptions, the fiber $\cFL_\gamma$ is a non-empty and locally finite dimensional ind-projective variety. 

There are two invariants one may attach to $\gamma$. For the first one, recall from \cite[\S1 Lemma 2]{KL:Fl} that the set of $G\laurent{t}$-conjugacy class of maximal tori in $G\laurent{t}$ is parameterized by the set $\ubar W$ of conjugacy classes in the finite Weyl group $W$. The \textit{type of $\gamma$} is the conjugacy class $[w] \in \ubar W$ that corresponds to the maximal torus $Z_{G\laurent{t}}(\gamma)$. For the second invariant, recall that the fiber $\cFL_\gamma$ consists of cosets $gI$ so that $\Ad(g)\inv \gamma \in \fI^\circ$. Since $\fI^\circ \subset \fg\series{t}$, by setting $t = 0$ we obtain an element $(\Ad(g)\inv\gamma)|_{t=0} \in \fg$ which is necessarily nilpotent. Let $\RT(\gamma) \subset \ubar \cN$ denote the set of nilpotent orbits that arise in this way, and let $\RTmin(\gamma)$ denote the set of minimal elements (with respect to the closure order) in $\RT(\gamma)$. The set $\RTmin(\gamma)$ is called the \textit{minimal reduction types of $\gamma$} \cite{Yun:RTmin}. It plays an important role in understanding affine Springer representations, see \cite{Jakob-Yun,Chua,Chua:RTmin} for examples.

In \cite{Lusztig:Phi,Lusztig:Phi2} Lusztig constructed a pair of maps
\begin{equation*}
	\Phi: \ubar W \surjects \ubar \cN,\quad \Psi: \ubar \cN \injects \ubar W
\end{equation*}
using the Bruhat decomposition of $G$, and $\Psi$ is a section of $\Phi$. Yun \cite{Yun:RTmin} showed that the maps $\Phi$, $\Psi$ are related to two maps
\begin{equation*}
	\RTmin: \ubar W \surjects \ubar \cN,\quad \KL: \ubar \cN \injects \ubar W
\end{equation*}
defined in a completely different way using rstn elements. More precisely, he showed that for any $[w] \in \ubar W$, there is a unique $\BO \in \ubar \cN$ such that $\BO \in \RTmin(\gamma)$ for any generic rstn element $\gamma$ of type $[w]$. Then one defines $\RTmin([w]) := \BO$. On the other hand, given any $\BO \in \ubar \cN$, Kazhdan-Lusztig showed in \cite{KL:Fl} that there is a unique $[w] \in \ubar W$ so that any generic element in $\BO + t \fg\series{t} \subset \fg\laurent{t}$ is rstn of type $[w]$. One sets $\KL(\BO) := [w]$. 

\begin{theorem}[{\cite[1.11 Theorem]{Yun:RTmin}}]
	The pair of maps $(\RTmin,\KL)$ coincides with $(\Phi,\Psi)$.
\end{theorem}

Define the \textbf{(Weyl group) cyclotomic level map} to be the composition
\begin{equation*}
	\cl_W := \cl_n \circ \RTmin: \ubar W \aro \BZ_{\ge 1},
\end{equation*}
so that we have a commutative triangle
\begin{equation*}
	\begin{tikzcd}
		\ubar W \ar[rr, "\cl_W"] \ar[dr, "\RTmin"']
		&& \BZ_{\ge 1}\\
		& \ubar \cN \ar[ur, "\cl_n"']
	\end{tikzcd}.
\end{equation*}

\begin{theorem}[{\cite[Theorem 3.3.5]{SYZ:cl}}]
	For any $w \in W$, $\cl_W([w])$ equals the maximal order of eigenvalues of $w$ on the reflection representation $W \acts \fh$. 
\end{theorem}

Alternatively, let $\operatorname{char}(w)$ denote the characteristic polynomial of $w \acts \fh$. Since the representation $W \acts \fh$ is defined over $\BQ$ and $w$ has finite order, $\operatorname{char}(w)$ is a product of cyclotomic polynomials $\Phi_n$. Then $\cl_W([w])$ equals the largest $n$ so that $\Phi_n$ appears as a factor of $\operatorname{char}(w)$. This explains the terminology. In fact, this was our original definition of cyclotomic level.

The theorem is equivalent to the following statement: if $\cl_W$ is defined as in the theorem, $\BO \in \ubar \cN$ is any orbit, and $[w]$ is any preimage of $\BO$ under $\RTmin$, then $\cl_W([w]) = \cl_n(\BO)$. In fact, any such $[w]$ can be realized as the saturation of a class $[w]_L$ in a Levi $L$, and if $L$ is chosen to be minimal, then $[w]_L$ is the image of an $\BO_L \subseteq \BO \cap \fl$ under the Levi version of the map $\KL$. Thus it suffices to show $\cl_W \circ \KL = \cl_n$. The proof of this is a modification of the argument in \cite[\S 9]{KL:Fl} where the statement is proven for distinguished $\BO$. Roughly speaking, the number $\cl_n(\BO)$ controls the minimal valuation of a generic lift $\gamma \in \BO + t \fg\series{t}$, and the latter contributes to the value $\cl_W([w])$.

\section{Main conjectures}\label{sec:conj}

We now turn to our conjectures. Let $\fg$ be of simply-laced type. Let $k \in \BZ$ be such that $k + \check \Bh := m \ge 1$. Let $\check \BO(m)$ be the orbit defined in Theorem \ref{thm:Om} in $\check \fg$. We expect the simple affine vertex algebra $L_k(\fg)$ plays a role in the affine setting similar to the primitive ideal quotient $\cU(\fg)/J_{\lambda,max}$ in the classical situation. In particular, the associated variety of $L_k(\fg)$ should be related to $\bd \check \BO(m)$, and the simple modules of $L_k(\fg)$ should be related to cells containing $w_m$. This cannot be literally true as there are known examples where $X_{L_k}$ is not contained in the nilpotent cone. Partially motivated by the constructions in \cite{AFK}, we believe the center of the Bala-Carter Levi of $\check \BO(m)$ should be responsible for the semisimple part in $X_{L_k}$.

\subsection{Conjecture on associated variety}

Write $\check \BO(m) = \Sat_{\check L}^{\check G} \BO_{\check L} = \Ad G\cdot \BO_{\check L}$ where $\check L$ is a Bala-Carter Levi of $\check \BO(m)$ (so that $\BO_{\check L}$ is distinguished in $\check \fl$). Let $\fp = \fl \oplus \fu$ be a parabolic subalgebra in $\fg$ containing the Langlands dual $\fl$ of $\check \fl$ as its Levi and $\fu$ as its radical, and let $\fz(\fl)$ denote the center of $\fl$. Write $\bd_L$ for the Barbasch-Vogan-Lusztig-Spaltenstein duality for the Levis $\check L$ and $L$.

\begin{conjecture}[{\cite[Conjecture 2.4.1]{SYZ:cl}}]\label{conj:AV}
	In the above setup, we have
	\begin{equation*}
		X_{L_k(\fg)} = \overline{\Ad G \cdot (\bd_L \BO_{\check L} \times \fz(\fl) \times \fu)}
	\end{equation*}
\end{conjecture}

All known cases are verified in \cite[\S 2.4]{SYZ:cl}. Note that the nilpotent part of the right hand side is 
\begin{align*}
	\overline{\Ad G \cdot (\bd_L \BO_{\check L} \times \fz(\fl) \times \fu)} \cap \cN 
	&= \overline{\Ad G \cdot (\bd_L \BO_{\check L} \times \fu)}\\
	&= \overline{\Ind_L^G \bd_L \BO_{\check L}}\\
	&= \overline{\bd \check \BO(m)}
\end{align*}
where $\Ind$ denotes Lusztig-Spaltenstein's induction of orbits \cite{Lusztig-Spaltenstein:induced}. Here the last equality follows from the equality $\bd \Sat_{\check L}^{\check G} = \Ind_L^G \bd_L$, see \cite[Proposition A2(c)]{Barbasch-Vogan:unipotent}.

Recall that $L_k$ is quasi-lisse if and only if the associated variety $X_{L_k}$ is contained in the nilpotent cone. From the above conjectural description, this happens if and only if $\fz(\fl) = 0$, if and only if $\check \fl = \check \fg$, if and only if $\check \BO(m)$ is distinguished. 

\begin{corollary}
	Suppose Conjecture \ref{conj:AV} holds. Then $L_k(\fg)$ is quasi-lisse if and only if $\check \BO(m)$ is distinguished.
\end{corollary}

\begin{example}
	In type $A_5$, the associated varieties predicted by the conjecture are summarized in the following table:
	\begin{center} 
		\begin{tabular}{ccccccc}
			\hline $m$ & $k$ & $\check \BO(m)$ & $\check L$ & $\BO_{\check L}$ & conjectural $X_{L_k}$ & confirmed by
			\\ \hline
			$\ge 6$ & $\ge0$ & $(6)$ & $A_5$ & reg & $\{0\}$ & well-known\\
			$5$ & $-1$  & $(51)$ & $A_4$ & reg & $\overline{\cS(A_4,\{0\})}$ & 
			\cite{Arakawa-Moreau:Omin,Jiang-Song}\\
			$4$ & $-2$ & $(42)$ & $A_3+A_1$ & reg & $\overline{\cS(A_3+A_1,\{0\})}$ & \\
			$3$ & $-3$ & $(33)$ & $2A_2$ & reg & $\overline{\cS(2A_2, \{0\})}$ & \cite{AFK}\\
			$2$ & $-4$ & $(222)$ & $3A_1$ & reg & $\overline{\cS(3A_1, \{0\})}$\\
			$1$ & $-5$  & $(1^6)$ & $T$ & reg & $\overline{\cS(T,\{0\})} =\fg$ & \cite{Gorelik-Kac}
		\end{tabular}
	\end{center}
\end{example}

\begin{example}
	In type $D_4$, the conjecture predicts the following:
	\begin{center}
		\begin{tabular}{ccccccc}
			\hline $m$ & $k$  & $\check \BO(m)$ & $\check L$ & $\BO_{\check L}$ & conjectural $X_{L_k}$ & confirmed by
			\\ \hline
			$\ge6$ & $\ge0$   & $(71)$ & $D_4$ & $(71)$ & $\{0\}$ & well-known\\
			$5$ & $-1$  & $(53)$ & $D_4$ & $(53)$ & $\overline{\BO_{min}}$ & \cite{Arakawa-Moreau:Omin}\\
			$4$ & $-2$   & $(53)$ & $D_4$ & $(53)$ & $\overline{\BO_{min}}$ & \cite{Arakawa-Moreau:Omin}\\
			$3$ & $-3$   & $(3311)$ & $A_2$ & $(3)$ & $\overline{\cS(A_2, \{0\})}$\\
			$2$ & $-4$  & $(3221)$ & $3A_1$ & $(2)^3$ & $\overline{\cS(3A_1, \{0\})}$\\
			$1$ & $-5$   & $(1^8)$ & $T$ & $\{0\}$ & $\overline{\cS(T,\{0\})} =\fg$ & \cite{Gorelik-Kac}
		\end{tabular}
	\end{center}
\end{example}

\subsection{Conjecture on simple modules}

We now turn to simple modules of $L_k(\fg)$. The details of this subsection and the next will appear in a future paper. Let $\xi_m \in \fh_{aff}^*$ and $w_m \in W_m \subset W_{aff}$ be the objects defined in \S\ref{subsec:cl-cells}. Recall from \S\ref{subsec:O(gaff)} that simple modules in the block $\cO_{k\Lambda_0}(\fg_{aff}) = \cO_{\xi_m - \hat \rho}(\fg_{aff})$ are parameterized by longest coset representatives of $W_{aff}/W_m$, or equivalently by those $y \in W_{aff}$ such that $y \le_L w_m$. In particular, we have an isomorphism of abelian groups
\begin{equation*}
	K_0 \cO_{\xi_m - \hat \rho}(\fg_{aff}) \cong \big( \cH_{aff}^\vee/ \cH_{aff, \not\le_L w_m}^\vee \big)|_{q=1}
\end{equation*}
where $\cH_{aff}$ denotes the affine Hecke algebra, and the notion of the dual module $\cH_{aff}^\vee$ and cells are defined as in \S \ref{subsec:KL-cells}. By using the affine version of twisting functors or by using Kashiwara-Tanisaki's localization, the above isomorphism can be made into $W_{aff}$-linear. Recall also that $\cO_{\xi_m - \hat \rho}(L_k(\fg))$ is a full subcategory of $\cO_{\xi_m - \hat \rho}(\fg_{aff})$. In particular, $K_0 \cO_{\xi_m - \hat \rho}(L_k(\fg))$ corresponds to a subspace in $\big( \cH_{aff}^\vee/ \cH_{aff, \not\le_L w_m}^\vee \big)|_{q=1}$.

\begin{conjecture}\label{conj:O(Lk)}
	Suppose $\check \BO(m)$ is distinguished. Then
	\begin{equation*}
		K_0 \cO_{\xi_m - \hat \rho}(L_k(\fg)) \cong \cH_{aff,\bc^L(w_m)}^\vee|_{q=1}.
	\end{equation*}
	In particular,
	\begin{equation*}
		\Irr \cO_{\xi_m - \hat \rho}(L_k(\fg)) \bijects \bc^L(w_m),\quad
		L(y \circ (\xi_m - \hat \rho)) \mathop{\mapsfrom} y. 
	\end{equation*}
\end{conjecture}

\begin{example}
	Let $\fg = \fso_8$, $k = -2$, $m = 4$. Then $w_m = s_2$, the reflection attached to the middle node in the extended Dynkin diagram. The affine left cell $\bc^L(s_2)$ consists of elements in $W_{aff}$ that have a unique reduced expression which ends with $s_2$ \cite[Chapter 12]{Bonnafe}. In other words,
	\begin{equation*}
		\bc^L(s_2) = \{s_2, s_0s_2, s_1s_2,s_3s_2, s_4s_2\}.
	\end{equation*}
	According to our conjecture, 
	\begin{equation*}
		\Irr \cO_{\xi_m - \hat \rho}(L_k(\fg)) = \{L(y \circ (\xi_m - \hat \rho)) \mid y \in \bc^L(s_2)\}.
	\end{equation*}
	This is proven by Perse \cite{Perse:D}.
\end{example}

\subsection{An affine picture}

Suppose $\check \BO(m)$ is distinguished. We will show in a future work that there exists an appropriate rstn element $\gamma \in \check \fg\laurent{t}$ which is a generic lift of an element of $\check \BO(m)$, such that there exists a $W_{aff}$-equivariant map from the right cell module $\cH_{aff,\bc^R(w_m)}|_{q=1}$ to the affine Springer representation $H^{top}(\check \cFL_\gamma)$ (here the $W_{aff}$-action on $H^{top}(\check \cFL_\gamma)$ is constructed by Lusztig \cite{Lusztig:aff-Sp-action}). Putting all the pieces together, we obtain the following affine analog to (\ref{diag:fin-duality}).

\begin{equation}\label{diag:aff-duality}
	\begin{tikzcd}
		\cH_{aff,\bc^R(w_m)}|_{q=1} \ar[d] \ar[rr, leftrightarrow, "\text{perfect pairing}", "(\ref{eqn:pairing})"']
		&& \cH_{aff,\bc^L(w_m)}^\vee|_{q=1} \ar[d, equal, "\text{Conjecture \ref{conj:O(Lk)}}"]\\
		H^{top}(\underbracket{\check \cFL_\gamma}) \ar[d, xshift=2ex, "\RTmin"']
		& \ubar \bc(w_m) \ar[dl, leftrightarrow, "\ref{thm:Om-wm}"]  \ar[ul, dotted] \ar[ur, dotted]
		& K_0 \cO_{\xi_m-\hat\rho}(\underbracket{L_k(\fg)}) \ar[d, xshift=2ex, "X_{(-)}"', "\text{Conjeture \ref{conj:AV}}"]\\		
		\check \BO(m) \ar[rr, leftrightarrow, "\bd"] 
		&& \bd \check \BO(m)
	\end{tikzcd}
\end{equation}

\printbibliography
\end{document}